\documentclass[twocolumn]{autart}    

\usepackage{tikz}
\usepackage{amsmath}
\usepackage{latexsym, amssymb, amsmath, amsbsy,amsopn, amstext}
\usepackage{graphicx,hyperref,booktabs}
\usepackage{cancel, color}
\usepackage{epstopdf,subcaption}
\graphicspath{{images/}}

\DeclareMathOperator{\Col}{Col}
\DeclareMathOperator{\Row}{Row}

\DeclareMathOperator{\argmax}{argmax}
\DeclareMathOperator{\BR}{BR}
\DeclareMathOperator{\SpanRow}{SpanRow}

\DeclareMathOperator{\diag}{diag}
\DeclareMathOperator{\lcm}{lcm}

\def\cal{\mathcal}

\def\ra{\rightarrow}

\def\a{\alpha}
\def\b{\beta}
\def\d{\delta}

\def\D{\Delta}

\def\ot{\otimes}

\def\M{\mbox{\boldmath$M$}}

\def\G{\Gamma}

\def\0{{\bf 0}}

\newcommand{\R}{{\mathbb R}}

\newcommand{\B}{{\mathcal B}}
\newcommand{\GM}{{\mathbb G}}

\def\dsum{\mathop{\sum}\limits}
\newtheorem{dfn}[thm]{Definition}
\newtheorem{prp}[thm]{Proposition}
\newtheorem{exa}[thm]{Example}

\begin{document}

\begin{frontmatter}

\title{Game Theoretic Control of Multi-Agent Systems\thanksref{footnoteinfo}}

\thanks[footnoteinfo]{ This work is supported partly by NNSF
    61333001 and 61273013 of China. Corresponding author: Daizhan Cheng. Tel.: +86 10 8254 1232.}

\author[1]{Ting Liu},
\author[2]{Jinhuan Wang},
\author[1]{Daizhan Cheng}


\address[1]{Laboratory of Systems and Control, AMSS,
Chinese Academy of Sciences, Beijing 100190, P.R.China}
\address[2]{School of Sciences, Hebei University of Technology,
Tianjin 300401, P.R.China}

\begin{keyword}
Finite game, (state based) potential function, multi-agent systems, designed utility function, semi-tensor product of matrices.
\end{keyword}

\begin{abstract}
Control of multi-agent systems via game theory is investigated. Assume a system level object is given, the utility functions for individual agents are designed to convert a multi-agent system into a potential game. First, for fixed topology, a necessary and sufficient condition is given to assure the existence of local information based utility functions. Then using local information the system can converge to a maximum point of the system object, which is a Nash equilibrium.
It is also proved that a networked evolutionary potential game is a special case of this multi-agent system. Second, for time-varying topology, the state based potential game is utilized to design the optimal control.  A strategy based Markov state transition process is proposed to assure the existence of state based potential function. As an extension of the fixed topology case, a necessary and sufficient condition for the existence of state depending utility functions using local information is also presented.  It is also proved that using better reply with inertia strategy, the system converges to a maximum strategy of the state based system object, which is called the recurrent state equilibrium.
\end{abstract}

\end{frontmatter}

\section{Introduction}

In recent years the game based control has received extensive attention from control community. The game theory has been applied to various control problems, including control of hybrid systems\cite{tom00}, planning hybrid power systems and distributed power control \cite{mei12}, \cite{hei06}, analysis and control of networks \cite{tem10}, distributed coverage of mobile agents \cite{yaz13}, \cite{zhu13}, road congestion control \cite{wan13}, just to name few.

\begin{figure}[!htbp]
\centering
\includegraphics[width=.45\linewidth]{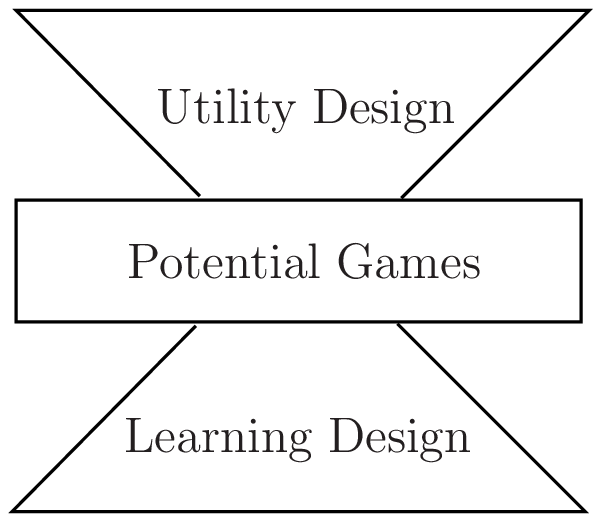}
\caption{Game Theoretic Approach \label{Fig.1}}
\end{figure}

Particularly, the game theoretic control is a promising new approach to the distributed control of multi-agent systems. In \cite{mar09} the consensus of multi-agent systems is investigated.
\cite{gop11} describes a solution framework for multi-agent control problem using game theory. Non-cooperative dynamic game theory provides an environment for formulating multi-agent decision control problems using distributed optimization \cite{bas99}. Learning is an important tool to realize a global goal for game theoretic approach of multi-agent systems \cite{abo14}.
An hourglass architecture was proposed in \cite{gop11} to illustrate the game theoretic control using potential games as the interface (refer to Fig. \ref{Fig.1}).

Consider a multi-agent system. Assume a system level objective function is given. This paper considers whether the system designer is able to design local information based utility functions for individual agents, such that the multi-agent system becomes a networked potential game with the objective function as the potential function. Then as the agents using their local information to maximize their utility functions the overall system can maximize the (system level) objective function. The paper consists mainly of two parts: First, the fixed topology case is considered. A necessary and sufficient condition is obtained to assure the existence of local information based utility functions. Then we show that the networked evolutionary potential games are special case of this kind of multi-agent systems.

Second, the case of time-varying topology is investigated. The concepts and results about state based potential game proposed by J.R. Marden \cite{mar12} have been used and combined with our previous results. Similar to fixed topology case, a necessary and sufficient condition is also obtained. Certain state transition process is proposed as strategy depending Markov process, which assures the non-decreasing requirement of the state based potential function. An illustrative example is presented, which shows that the game theoretic control makes the multi-agent system reach consensus with probability $1$.

The rest of this paper is organized as follows: Section 2 provides some necessary preliminaries, including (i) semi-tensor product of matrices, which is a basic tool in our approach; (ii)the potential equation, which is a key result used in this paper; (iii) networked potential game, which is an important example of the multi-agent systems concerned in this paper. Section 3 considers the game theoretic control for fixed topology multi-agent systems. A necessary and sufficient condition is presented. It is also shown that networked potential games meet the requirement. The game theoretic control for time-varying topology is discussed in Section 4. State based potential game is reviewed first. A necessary and sufficient conditions is also obtained. An illustrative example is presented. Section 5 is a brief conclusion.

\section{Preliminaries}

\subsection{Semi-tensor Product of Matrices}

This subsection gives a brief review for semi-tensor product (STP). We refer to \cite{che07}, \cite{che11}, \cite{che12} for details.
For statement ease, we first introduce some notations:

\begin{enumerate}
\item ${\cal M}_{m\times n}$: the set of $m\times n$ dimensional real matrices.
\item $\Col(A)$ ($\Row(A)$): the set of columns (rows) of ~$A$; $\Col_i(A)$ ($\Row_i(A)$): the $i$-th column (row) of ~$A$.
\item ${\cal D}_k=\{1,2,\cdots,k\}$,  ${\cal D}:={\cal D}_2$;
\item $\d_n^i$: the $i$-th column of the identity matrix ~$I_n$;
\item $\D_n=\{\d_n^i\,|\,i=1,2,\cdots,n\}$;
\item $L\in {\cal M}_{m\times r}$ is a logical matrix, if ~$\Col(L)\subset \D_m$. The set of  $m\times r$ logical matrices is denoted as ~${\cal L}_{m\times r}$;
\item Let ~$L$ be a logical matrix, i.e.,  ~$L\in {\cal L}_{m\times r}$. Then ~$L=[\d_m^{i_1}\;\d_m^{i_2}\;\cdots\; \d_m^{i_r}]$. For brevity,
$$
L=\d_m[i_1\;i_2\;\cdots\;i_r].
$$
\item  Set of random vectors:
$$
\varUpsilon_k:=\left\{(r_1,r_2,\cdots,r_k)^T\;\big|\; r_i\geq 0,~\dsum_{i=1}^kr_i=1 \right\}.
$$
\item  Set of random matrices:
$$
\varUpsilon_{m\times n}:=\left\{M\in {\cal M}_{m\times n}\;\big|\; \Col(M)\subset \varUpsilon_m \right\}.
$$
\end{enumerate}
\begin{dfn}\label{d2.1.1} Let $A\in {\cal M}_{m\times n}$, $B\in {\cal M}_{p\times q}$, and $t=\lcm(n,p)$ be
the least common multiple of $n$ and $p$. Then the (left) STP of $A$ and $B$, denoted by $A\ltimes B$, is
defined as
\begin{align}\label{2.1.1}
A\ltimes B:=\left(A\otimes I_{t/n}\right)\left(B\otimes I_{t/p}\right),
\end{align}
where the $\otimes$ is the Kronecker product.
\end{dfn}

\begin{rem}\label{r2.1.2}
\begin{enumerate}
\item If ~$n=p$, the STP defined in Definition \ref{d2.1.1} degenerates to the conventional matrix product. Hence, STP is a generalization of the conventional matrix product. Hence, we adopt as the convention that $AB:=A\ltimes B$.
\item All the major properties of conventional matrix product, such as the associativity and the distributivity etc., remain available.
\end{enumerate}
\end{rem}

The following proposition shows that the STP has certain communicative property.

\begin{prp}\label{p2.1.3} Given ~$A\in {\cal M}_{m\times n}$.
\begin{enumerate}

\item[{\rm 1.}] Let ~$Z\in \R^t$ be a column vector. Then
\begin{align}
\label{2.1.2} ZA=( I_t\otimes A)Z.
\end{align}

\item[{\rm 2.}] Let ~$Z\in \R^t$ be a row vector. Then
\begin{align}
\label{2.1.3} AZ=Z(I_t\otimes A).
\end{align}
\end{enumerate}
\end{prp}

To explore further communicating properties, we introduce the swap matrix.

\begin{dfn}\label{d2.1.4} A swap matrix ~$W_{[m,n]}\in {\cal M}_{mn\times
mn}$ is defined as follows:
\begin{align}
\label{2.1.4}
\begin{array}{ccr}
W_{[m,n]}&=&\d_{mn}[1,m+1,\cdots,(n-1)m+1;\\
~&~&2,m+2,\cdots,(n-1)m+2;\\
~&~&\cdots~;~m,2m,\cdots,nm].
\end{array}
\end{align}
\end{dfn}
The following proposition shows that the swap matrix is orthogonal.

\begin{prp}\label{p2.1.5}
\begin{align} \label{2.1.5} W_{[m, n]}^T=W_{[m, n]}^{-1}=W_{[n,m]}.
\end{align}
\end{prp}

Its fundamental function is to swap two factors.

\begin{prp}\label{p2.1.6}
\begin{enumerate}
\item[{\rm 1.}] Let ~$X\in \R^m$, $Y\in \R^n$ be two column vectors. Then
\begin{align} \label{2.11} W_{[m, n]}\ltimes X\ltimes Y=Y\ltimes X.
\end{align}
\item[{\rm 2.}] Let ~$X\in \R^m$, $Y\in \R^n$ be two row vectors. Then
\begin{align} \label{2.12} X\ltimes Y\ltimes W_{[m, n]}=Y\ltimes X.
\end{align}
\end{enumerate}
\end{prp}

\subsection{Potential Games}

\begin{dfn} \label{d2.2.1} \cite{gib92} A finite game is denoted by $G=\{N,S,c\}$, where (i) $N=\{1,\cdots,n\}$ is the set of players; (ii) $S=\prod_{i=1}^nS_i$ is the profile of strategies (or actions) $S_i$, where $S_i=\{1,\cdots,k_i\}$ is the set of strategies of player $i$; (iii) $c=(c_1,\cdots,c_n)$ and each $c_i:S\ra \R$ is the utility (or payoff) function of player $i$.
\end{dfn}

To use STP, the action $a_i\in S_i$ is denoted as $S_i=\left\{\d_{k_i}^j\;\big|\;j=1,\cdots,k_i\right\}$, where $j\in S_i$ is expressed as $j\sim \d_{k_i}^j$.
Then each utility function $c_i$ can be expressed as
\begin{align}\label{2.2.1}
c_i=V^c_i\ltimes_{j=1}^na_j,\quad i=1,\cdots,n,
\end{align}
where $V^c_i\in \R^{k}$ is called the structure vector of $c_i$ ($k:=\prod_{i=1}^nk_i$).

The set of finite games with $|N|=n$, $|S_i|=k_i$ is denoted by ${\cal G}_{[n;k_1,\cdots,k_n]}$. $G\in {\cal G}_{[n;k_1,\cdots,k_n]}$ is completely determined by $V_G:=(V^c_1,\cdots,V^c_n)\in \R^{nk}$. Hence, ${\cal G}_{[n;k_1,\cdots,k_n]}$ has a vector space structure as $\R^{nk}$ \cite{chepr}.

\begin{dfn} \label{d2.2.2} Consider a finite game $G=\{N,S,C\}$.
$G$ is a potential game if there exists a function $P:S\ra \R$, called the potential function, such that for every $i\in N$ and for every $s_{-i}\in S_{-i}:=\prod_{j\neq i}S_j$ and $\forall \a,\b\in S_i$
\begin{align}\label{2.2.2}
c_i(\a, s_{-i})-c_i(\b,s_{-i})=P(\a, s_{-i})-P(\b,s_{-i}).
\end{align}
\end{dfn}

Potential games have some nice properties, which make them helpful in control design. We listed some as follows:

\begin{thm}\label{t2.2.3}  \cite{mon96}
If $G$ is a potential game, then the potential function $P$ is unique up to a constant number. Precisely, if $P_1$ and $P_2$ are two potential functions of $G$, then $P_1-P_2=c_0\in \R$.
\end{thm}

\begin{thm} \label{t2.2.4}  \cite{mon96}
Every finite potential game possesses a pure Nash equilibrium. The  myopic best response adjustment (MBRA) leads to a Nash equilibrium. (Please also refer to (\ref{2.3.4})-(\ref{2.3.5}) for MBRA.)
\end{thm}

Consider a finite game $
G=\{N,S,c\}$, where $|N|=n$, $|S_i|=k_i$, $i=1,\cdots,n$, and the payoff function of player $i$ is denoted as (\ref{2.2.1}).
We need some new notations as:
\begin{align}\label{2.2.3}
k^{[p,q]}:=
\begin{cases}
\prod_{j=p}^qk_j,& q\geq p\\
1,               & q<p;
\end{cases}
\end{align}
and
\begin{align}\label{2.2.4}
E_i:=I_{k^{[1,i-1]}}\otimes {\bf 1}_{k_i}\otimes I_{k^{[i+1,n]}}\in {\cal M}_{k\times k/k_i},~ i=1,\cdots,n.
\end{align}

Then we construct a linear equation, called the potential equation, as
\begin{align}\label{2.2.5}
\begin{bmatrix}
-E_1&E_2&0&\cdots&0\\
-E_1&0&E_3&\cdots&0\\
\vdots&~&~&~&~\\
-E_1&0&0&\cdots&E_n
\end{bmatrix}\begin{bmatrix}
\xi_1\\
\xi_2\\
\vdots\\
\xi_n
\end{bmatrix}=
\begin{bmatrix}
(V^c_2-V^c_1)^T\\
(V^c_3-V^c_1)^T\\
\vdots\\
(V^c_n-V^c_1)^T\\
\end{bmatrix},
\end{align}
where $\xi_i\in \R^{k/k_i}$. Then we have the following result:

\begin{thm} \cite{che14}, \cite{chepr} \label{t2.2.5} A finite game $G\in {\cal G}_{[n;k_1,\cdots,k_n]}$ is a potential game, if and only if the potential equation (\ref{2.2.5}) has solution. Moreover, if $\xi$ is a solution then the potential function can be expressed as
\begin{align}\label{2.2.6}
P(x_1,\cdots,x_n)=V^P\ltimes_{j=1}^nx_j,
\end{align}
where $V^P$, the structure vector of the potential function, is
\begin{align}\label{2.2.7}
V^P=V^c_1-\xi^T_1E_1^T.
\end{align}
\end{thm}

\subsection{Networked Evolutionary Games}

\begin{dfn}\label{d2.3.1}\cite{che15} A networked evolutionary game (NEG), denoted by $\Gamma=\{(N,E), G, \Pi\}$, consists of
\begin{enumerate}
\item[(i)] a network graph  $(N, E)$;
\item[(ii)] a fundamental network game (FNG), $G$, such that if $(i,j)\in E$, then $i$ and $j$ play FNG with strategies $a_i(t)$ and $a_j(t)$ respectively;
\item[(iii)] a local information based strategy updating rule (SUR), $\Pi$.
\end{enumerate}
\end{dfn}

A Markov-type strategy profile dynamics of an NEG can be expressed as
\begin{align}\label{2.3.1}
\begin{cases}
a_1(t+1)=f_1(a_1(t),\cdots,a_n(t))\\
a_2(t+1)=f_2(a_1(t),\cdots,a_n(t))\\
\vdots\\
a_n(t+1)=f_n(a_1(t),\cdots,a_n(t)),
\end{cases}
\end{align}
where $a_i(t)\in {\cal D}_{k_i}$, $i=1,\cdots,n$.

Using vector form expression of the strategies, (\ref{2.3.1}) can be expressed into its algebraic state space form as
\begin{align}\label{2.3.2}
a(t+1)=La(t),
\end{align}
where $a(t)=\ltimes_{i=1}^na_i(t)$, and when only pure strategies are allowed then
$L\in {\cal L}_{k\times k}$; if the mixed strategies are allowed then $L\in \varUpsilon_{k\times k}$. ($k=\prod_{i=1}^nk_i$.)

The $f_i$ in (\ref{2.3.1}) is determined by the SUR. In the following, MBRA, as the only SUR used in this paper, is briefly introduced.  Denote by $U(i)$ the neighborhood of node $i$. As a convention, we assume $i\in U(i)$.

Let $c_{ij}$ be the payoff of the player $i$ in the FNG between $i$ and $j$, then the overall payoff of player $i$ is
\begin{align}\label{2.3.3}
c_i(t)=\dsum_{j\in U(i)\backslash \{i\}}c_{ij}, \; i\in N.
\end{align}

Assume the SUR used is MBRA, that is
\begin{align}\label{2.3.4}
a_i(t+1)\in \argmax_{a_i\in S_i} (c_i(a_i,a_{-i}(t)):=\BR_i(t),
\end{align}
where the elements of $BR_i(t)$ are named as
$$
\BR_i(t):=\left\{a^i_1,a^i_2,\cdots,a^i_{\mu_i(t)}\right\}.
$$
$a_i(t+1)$ is precisely defined as follows:
\begin{itemize}
\item If $a_i(t)\in \BR_i(t)$, then
$$
a_i(t+1)=a_i(t);
$$
\item If $a_i(t)\not\in \BR_i(t)$, then
$$
\begin{array}{l}
a_i(t+1)=a^i_j, ~\mbox{with probability}~ \frac{1}{\mu_i(t)},\\
~~~~~~~~~~~~~~~~~~~~~~~~~~~~ j=1,\cdots,\mu_i(t).
\end{array}
$$
\end{itemize}

If only the local information is allowed to use, then (\ref{2.3.2}) becomes

\begin{align}\label{2.3.5}
a_i(t+1)\in \argmax_{a_i\in U(i)} c_i(a_i,a_{U(i)\backslash \{i\}}(t)),
\end{align}
and the rest determining process remains unchanged.

\begin{rem}\label{r2.3.2} If only the local information is allowed to use, the convergence described in Theorem \ref{t2.2.4} is not assured, unless the utilities are themselves local information depending. That is why we need to design the local information depending utilities.
\end{rem}

\section{Utility Design}

As pointed in \cite{gop11} or \cite{mar12}, in general, a multi-agent system may not originally a networked evolutionary game. There is a system level objective function $\phi:S\ra \R$ that a system designer seeks to maximize. Then we may design a set of suitable utility functions such that the system becomes a potential game with $\phi$ as its potential function. This is the first task in game theoretic approach, which is described in top part of Fig. \ref{Fig.1}. Second task is: since the information used is a local one, a learning SUR is necessarily to assure the system converges to the maximum object value. We start with a motivated example.

\subsection{A Motivated Example}

\begin{figure}[!htbp]
\centering
\includegraphics[width=.45\linewidth]{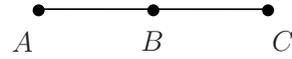}
\caption{A Network\label{Fig.2}}
\end{figure}

 Assume there are three players: $A,~B,~C$, connected as in Fig. \ref{Fig.2}. Assume each player has two strategies $S_i=\{1,2\}$, $i=1,2,3$. The system level object is to guide players to a synchronized form $a_1(t)=a_2(t)=a_3(t)=1$. So a system level objective function can be described as
$$
\phi(t):=\left|\{i\;\big|\;a_i(t)=1\}\right|.
$$
The system designer wants to maximize $\phi$, which leads the system to the synchronization.

The structure vector of $\phi$ is figured out as
$$
V^{\phi}=[3,2,2,1,2,1,1,0].
$$
Assume
\begin{align}\label{3.1.1}
\begin{array}{ccl}
c_1(a)&=&[2,1,1,0]a_1a_2\\
c_2(a)&=&[3,4,2,3,2,0,1,-1]a_1a_2a_3\\
c_3(a)&=&[1,0,1,0]a_2a_3.
\end{array}
\end{align}
It is easy to calculate that
$$
\begin{array}{ccl}
V^c_1&=&[2,2,1,1,1,1,0,0]\\
V^c_2&=&[3,4,2,3,2,0,1,-1]\\
V^c_3&=&[1,0,1,0,1,0,1,0].
\end{array}
$$
Using (\ref{2.2.4}), we have
$$
\begin{array}{ccl}
E_1&=&{\bf 1}_2\otimes I_4\\
E_2&=&I_2\otimes {\bf 1}_2\otimes I_2\\
E_3&=&I_4\otimes {\bf 1}_2.
\end{array}
$$
Then the potential equation (\ref{2.2.5}) can be built, and a solution is obtained as
$$
\xi=[0,1,0,1,1,3,1,0,-1,0,0,1]^T,
$$
and $\xi_1$ follows as
$$
\xi_1=[0,1,0,1]^T.
$$
According to Theorem \ref{t2.2.5} this game is potential.
Moreover, using (\ref{2.2.7}) we can calculate that
$$
V^P=V^c_1-\xi_1^TE_1^T=[2,1,1,0,1,0,0,-1].
$$
Note that
$$
P(a)=V^Pa_1a_2a_3=\phi(a)-1,
$$
hence $\phi(a)$ is also a potential function.

Observing that in this networked game the utility functions $c_i$ depends only on its neighborhood. Hence the local information MBRA, as defined in (\ref{2.3.5}) is the same as the global information MBRA, as in (\ref{2.3.4}). According to Theorem \ref{t2.2.4}, the local information MBRA can lead the system to a Nash equilibrium, which maximizes the system objective function.

Motivated by this example, it is natural to seek a set of local information based utility functions, which then can lead the system to a maximum point.

\subsection{Local Information Based Utility Functions}

Consider a game $G\in {\cal G}_{[n;k_1,\cdots,k_n]}$. Assume $U\subset N$, say $U$ could be a neighborhood of a node. Then we try to ``draw" the nodes of $U$ from all the nodes in $N$. We construct a matrix, called the $U$-drawing matrix, to do this. Set
\begin{align}\label{3.2.1}
\Gamma_U:=\otimes_{i=1}^n\gamma_i,
\end{align}
where
\begin{align}\label{3.2.2}
\gamma_i=
\begin{cases}
I_{k_i},\quad i\in U\\
{\bf 1}_{k_i}^T,\quad \mbox{Otherwise}.
\end{cases}
\end{align}
Then we have the following result.

\begin{lem}\label{l3.2.1} Let $U\subset N$. Then
\begin{align}\label{3.2.3}
\ltimes_{j\in U}a_j=\Gamma_U\ltimes_{i=1}^na_i.
\end{align}
\end{lem}

\begin{rem}\label{r3.2.2}
\begin{enumerate}
\item The proof of Lemma \ref{l3.2.1} is based on the fact that for two column vectors $X,~Y$,
$$
X\ltimes Y=X\otimes Y.
$$
Then a straightforward computation leads to the conclusion.

\item From the proof it is clear that formula (\ref{3.2.3}) is also available for $a_i\in \varUpsilon_{k_i}$, $i=1,\cdots,n$.
\end{enumerate}
\end{rem}

Consider $G\in {\cal G}_{[n;k_1,\cdots,k_n]}$ being a networked evolutionary game. Assume the utilities are adjustable. A natural question is: are we able to design a neighborhood-determinant utilities such that $G$ becomes a potential game and local information is enough to assure the convergence.

\begin{thm}\label{3.2.3} Assume  $G\in {\cal G}_{[n;k_1,\cdots,k_n]}$ is a utility-adjustable networked evolutionary game. The system objective function is
$$
\phi(a)=V^{\phi}\ltimes_{i=1}^na_i.
$$
Then there exists a set of neighborhood-determinant utilities, which turn $G$ to be a potential game, if and only if
\begin{align}\label{3.2.4}
V^{\phi}\in \bigcap_{i=1}^n \SpanRow \begin{bmatrix}
\Gamma_{U(i)}\\
E_i^T
\end{bmatrix}.
\end{align}
\end{thm}

\noindent{\it Proof}. By the requirement, we have
\begin{align}\label{3.2.5}
c_i(a)=V^c_i\ltimes_{j\in U(i)}a_j=V^c_i\Gamma_{U(i)}a,\quad i=1,\cdots,n,
\end{align}
where $a=\ltimes_{i=1}^na_i$.

Using the same argument as in the proof of Theorem \ref{t2.2.5}, one sees that $G$ is a potential game, if and only if,
\begin{align}\label{3.2.6}
\phi(a)-c_i(a)=d_i(a_{-i}),\quad i=1,\cdots,n,
\end{align}
where $a_{-i}\in S_{-i}$, which means $d_i$ is independent of $a_i$.

Plugging (\ref{3.2.5}) into (\ref{3.2.6}) yields
\begin{align}\label{3.2.7}
V^{\phi}-V^c_i\Gamma_{U(i)} =V^d_i E_i^T,\quad i=1,\cdots,n,
\end{align}
which is equivalent to (\ref{3.2.4}).
\hfill $\Box$

For convenience in use, denote by $B$ a basis of the vector space on the right hand side of (\ref{3.2.4}). That is, $\Row(B)$ are linearly independent and
\begin{align}\label{3.2.8}
\SpanRow (B)=
\bigcap_{i=1}^n \SpanRow \begin{bmatrix}
\Gamma_{U(i)}\\
E_i^T
\end{bmatrix}.
\end{align}

\subsection{Networked Potential Games}

Consider an evolutionary game $\GM=\{(N,E), G, \Pi\}$, assume the fundamental network game $G$ is potential, then we call $\GM$ a networked potential game (NPG). Assume the system objective function
$$
\phi(a)=\dsum_{e\in E}P_e(a),
$$
where $P_e$ is the potential function for $G$ over the edge $e$. Using the natural utility function (\ref{2.3.3}),  it is clear that (\ref{3.2.4}) is satisfied. This fact shows that pretty of multi-agent systems verify the requirements of Theorem \ref{t3.2.3}.

In the following we give an example.

\begin{exa}\label{e3.3.1}
Consider an NEG $\GM=\left\{\left(N, E\right), G, \Pi \right\}$, where the network graph is  as in Fig. \ref{fig3}, and the FNG, $G$, is the Prisoner's Dilemma game with payoff matrix shown in Table \ref{tab2}.
\begin{figure}[!htbp]
\centering
\includegraphics[width=.35\linewidth]{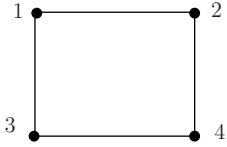}
\caption{Network Graph, $S_4$\label{fig3}}
\end{figure}

\begin{table}[h]
  \centering
  \caption{Payoff Bi-matrix of Prisoner's Dilemma \label{tab2}}
  \begin{tabular}{ccc}
    \toprule $c_1\backslash c_2$ & $1(cooperate)$ & $2(defect)$\\
    \midrule $1(cooperate)$ & 3,~3 & 0,~5\\
    $2(defect)$ & 5,~0 & 1,~1\\
    \bottomrule
  \end{tabular}
\end{table}
Using potential equation, we can prove that Prisoner's Dilemma is a potential game and one of its potential function is $P(a_1,a_2) =[-2,0,0,1]a_1a_2$. Then the overall potential function $\phi$ is
\begin{align*}
\begin{array}{l}
\phi(a)\\
=P(a_1, a_2)+P(a_1, a_3)+P(a_2, a_4)+P(a_3, s_4)\\
=V^PI_4\ot {\bf 1}_4^Ta +V^P I_2\ot {\bf 1}_2^T\ot I_2\ot {\bf 1}_2^T a \\
\qquad+ V^P {\bf 1}_2^T\ot I_2\ot{\bf 1}_2^T\ot I_2 a +V^P {\bf 1}_4^T\ot I_4 a,\\
:=V^{\phi} a,
\end{array}
\end{align*}
where
\[
V^{\phi}=[ -8  , -4  , -4  , -1  , -4  , 0  , -1  , 2  , -4  , -1  , 0  , 2  , -1  , 2  , 2  , 4].
\]
According to (\ref{3.2.1}) and (\ref{3.2.2}), we have
\begin{align*}
\begin{array}{ll}
\G_{U(1)}= I_8\ot {\bf 1}_2^T;& \G_{U(2)}= I_4\ot {\bf 1}_2^T\ot I_2; \\
\G_{U(3)}= I_2\ot {\bf 1}_2^T\ot I_4;  &\G_{U(4)}= {\bf 1}_2^T\ot I_8.
\end{array}
\end{align*}
Additionally,
\begin{align*}
\begin{array}{ll}
E_1^T= {\bf 1}_2^T\ot I_8;& E_2^T= I_2\ot {\bf 1}_2^T\ot I_4; \\
E_3^T= I_4\ot {\bf 1}_2^T\ot I_2;  &E_4^T= I_8\ot {\bf 1}_2^T.
\end{array}
\end{align*}
Define
\[
\G=\begin{bmatrix}
\G_{U(1)}&~&~&\\
&\G_{U(2)}&~&~\\
 &~&\G_{U(3)}&\\
  &~&~&\G_{U(4)})\\
\end{bmatrix},
\B =\begin{bmatrix}
E_1^T&~&~&\\
&E_2^T&~&~\\
 &~&E_3^T&\\
  &~&~&E_4^T\\
\end{bmatrix}.
\]
It is easy to check that
\[
{\bf 1}^T_4\ot V^{\phi}\in \SpanRow \begin{bmatrix}\G\\ \B\end{bmatrix},
\]
which is equivalent to (\ref{3.2.4}).
\end{exa}

\section{State Based Potential Games}

\subsection{A Brief Review}

The state based potential game is proposed in \cite{mar12}. This subsection briefly review some related basic concepts and results.

\begin{dfn}\label{d4.1.1}
\begin{enumerate}
\item A finite state based evolutionary game is a tuple $G=\{N,S,C,X,P\}$, where $N=\{1,2,\cdots,n\}$ is the set of players; $S=\prod_{i=1}^nS_i$ is the strategy profile; $C=\{c_1,\cdots,c_n\}$ is the set of payoff (utility) functions, and $c_i:X\times S\ra \R$; $X$ is the state space; $P: X \times S\ra  \triangle(X)$ is the the probability distributions over the finite state space $X$.
\item $X=\{x_1,\cdots,x_r\}$, where $r=\left|X\right|$.  Similar to $a_i\in \D_{k_i}$, we can also to express states in vector form as
$$x_i\sim \d_r^i,\quad i=1,\cdots,r.
$$
\item The state $x(t)$ satisfies
\begin{align}\label{4.1.3}
x(t+1)\sim P\left(x(t),a(t)\right),
\end{align}
which can be expressed into its algebraic state space form as
\begin{align}\label{4.1.301}
x(t+1)=M_Px(t)a(t),
\end{align}
where $a(t)=\ltimes_{i=1}^na_i(t)$.
\item The strategy dynamics of an evolutionary game is of the following form
\begin{align}\label{4.1.1}
\begin{cases}
a_{1}(t+1)=f_1(x(t+1),a(t),c(t))\\
a_{2}(t+1)=f_2(x(t+1),a(t),c(t))\\
\vdots\\
a_{n}(t+1)=f_n(x(t+1),a(t),c(t)),
\end{cases}
\end{align}
where $c(t)=(c_1(t),c_2(t),\cdots,c_n(t))$.
The dynamics is determined by an SUR. In this paper we assume $f_i$ is independent of $c$, hence it can be expressed in algebraic state space form as
\begin{align}\label{4.1.2}
a(t+1)=M_Fx(t+1)a(t).
\end{align}

\end{enumerate}
\end{dfn}

\begin{dfn}\label{d4.1.2} The action state pair $[a^*,x^*]$ is a recurrent state equilibrium with respect to the state transition process $P(\cdot)$ if the following two conditions are satisfied:
\begin{itemize}
\item The state $x^*$ satisfies $x^*\in X(a^*|x)$ for every state $x\in X(a^*|x^*)$.
\item For every agent $i\in N$ and every state $x\in X(a^*|x^*)$,
$$c_i(x, a_i^*,a^*_{-i})\geq c_i(x, a_i,a^*_{-i}),\ \forall a_i\in S_i$$
\end{itemize}
\end{dfn}

\begin{dfn}\label{d4.1.3} A state based game $G=\{N,S,C,X,P\}$ is a state based potential game if there exists a potential function $\phi:X\times S\ra \R$ that satisfies the following two properties for every action state pair $[a,x]\in S\times X$:
\begin{itemize}
\item For any agent $i\in N$ and action $a_i'\in S_i$
\begin{align}\label{4.1.4}
c_i(x, a_i',a_{-i})-c_i(x, a)=\phi(x, a_i',a_{-i})-\phi(x,a).
\end{align}
\item For any state $x'$ in the support of $P(x,a)$
\begin{align}\label{4.1.5}
\phi(x',a)\geq \phi(x,a).
\end{align}
\end{itemize}
\end{dfn}

The ``better reply with inertia" dynamics is important for state based games. Define an agent's strict better reply set for any action state pair $[a,x]\in S\times X$ as
$$
B_i(x,a):=\left\{a_i'\in S_i\;\big|\; c_i(x,a_i',a_{-i})>c_i(x,a)\right\}.
$$
The better reply with inertia dynamics can be described as follows.
\begin{itemize}
\item If $B_i(x(t), a(t-1))=\emptyset$ then
\begin{align}\label{4.1.6}
p^{a_i(t-1)}_i=1.
\end{align}
\item Otherwise, if $B_i(x(t), a(t-1))\neq \emptyset$ then
\begin{align}\label{4.1.7}
\begin{cases}
p^{a_i}_i=\epsilon,\quad a_i=a_i(t-1)\\
p^{a'_i}_i=\frac{(1-\epsilon)}{\left|B_i(x(t), a(t-1))\right|},\quad a'_i\in B_i(x(t), a(t-1))\\
p_i^{a''_i}=0,\quad \mbox{Otherwise},
\end{cases}
\end{align}
where $\epsilon \in (0,1)$ is referred to as the agent's inertia.
\end{itemize}

\begin{thm}\label{t4.1.4} Let $G=\left\{N,S,C,X,P\right\}$ be a state based potential game with potential function $\phi:X\times S\ra \R$. If all agents adhere to the better reply with inertia dynamics then the action state pair converges almost surely to an action invariant set of recurrent state equilibria.
\end{thm}

\subsection{State-depending Utility Design}

Assume a state based multi-agent system is described as $\G=\left\{N,S,c,X,P\right\}$, where $N=\{1,\cdots,n\}$; $S=\prod_{i=1}^nS_i$ and $S_i=\{1,\cdots,k_i\}$, $i=1,\cdots,n$; $c_i:X\times S\ra \R$; $P:X\times S\ra \overline{X}$; $X=\{x_1,\cdots,x_r\}$; and  the system level object is a state-depending function $\phi=\phi(x,a)$.

The system objective function can be expressed as
\begin{align}\label{4.2.1}
\phi(x,a)=V^{\phi}xa,
\end{align}
where $a=\ltimes_{i=1}^na_i$, $V^{\phi}\in \R^{rk}$, $k=\prod_{i=1}^nk_i$.

Split $V^{\phi}$ into $r$ equal blocks as
$$
V^{\phi}=[V^{\phi}_1,V^{\phi}_2,\cdots,V^{\phi}_r].
$$
Assume $x=x_i$ is fixed. Then the structure vector of $\phi(x_i,\cdot)$ is
$$
V^{\phi(x_i,\cdot)}=V^{\phi}_i,\quad 1\leq i\leq r.
$$
Corresponding to  $\phi(x_i,\cdot)$, we can construct $B_i$ as in (\ref{3.2.8}), that is
$$
\SpanRow (B_i)=
\bigcap_{j=1}^n \SpanRow \begin{bmatrix}
\Gamma_{U^i(j)}\\
E_j^T
\end{bmatrix},
$$
where $U^i(j)$ is the neighborhood of player $j$ under fixed state $x=x_i$.
Define
\begin{align}\label{4.2.1}
B:=\diag(B_1,B_2,\cdots,B_r).
\end{align}

Next, we need to design the state evolutionary process (SEP). We suggest the following two ways to construct it.

\begin{itemize}
\item SEP-1 (Remaining Priority):

Construct
$$
BR(x(t),a(t)):=\left\{x_j\;\big|\; \phi(x_j,a(t))>\phi(x(t),a(t))\right\}.
$$
Then
\begin{itemize}
\item if $BR(x(t),a(t))=\emptyset$, then
\begin{align}\label{4.2.2}
x(t+1)=x(t);
\end{align}
\item if $BR(x(t),a(t))\neq\emptyset$, then
\begin{align}\label{4.2.3}
\begin{array}{l}
x(t+1)=x_j,\quad \mbox{with probability } p_j=\frac{1}{\left| BR(x(t),a(t))   \right|}\\
~~~~~~~~~~~~~~~~~~~~~~~~~~~~~~~~~x_j\in BR(x(t),a(t)).
\end{array}
\end{align}
\end{itemize}

\item SEP-2 (Equal Probability):
$$
BR(x(t),a(t)):=\left\{x_j\;\big|\; \phi(x_j,a(t))\geq \phi(x(t),a(t))\right\}.
$$
Then
\begin{align}\label{4.2.301}
\begin{array}{l}
x(t+1)=x_j, \text{with probability } p_j=\frac{1}{\left| BR(x(t),a(t))   \right|}\\
~~~~~~~~~~~~~~~~~~~~~~~~~~~~~~~~~x_j\in BR(x(t),a(t)).
\end{array}
\end{align}
\end{itemize}

By the above construction it is clear that

\begin{lem}\label{l4.2.1}
Both the SEP-1 determined by (\ref{4.2.2})-(\ref{4.2.3}) and the SEP-2 determined by (\ref{4.2.301})  assure (\ref{4.1.5}).
\end{lem}

Similar to Theorem \ref{3.2.3}, we have the following result:

\begin{thm}\label{t4.2.2}  Assume  $\Gamma=\left\{N,S,c,X,P\right\}$  is a state based utility-adjustable networked evolutionary game. The system objective function is
$$
\phi(x, a)=V^{\phi}x\ltimes_{i=1}^na_i.
$$
Then there exists a set of neighborhood-determinant state based utilities, which turn $\Gamma$ to be a state based potential game, if and only if
\begin{align}\label{4.2.4}
V^{\phi}\in \SpanRow B.
\end{align}
\end{thm}

\begin{rem}\label{r4.2.3} Using the SEP-1 determined by (\ref{4.2.2})-(\ref{4.2.3}) or the SEP-2 determined by (\ref{4.2.301}) and the system objective function satisfying (\ref{4.2.4}), we can convert a multi-agent system to a state based potential game. Then using better reply with inertia dynamics, the system will converge almost surely to an action invariant set of recurrent state equilibria.
\end{rem}

\subsection{An Illustrative Example}

\begin{exa}\label{e4.3.1}

\begin{figure}[!htbp]
\centering
\includegraphics[width=.45\linewidth]{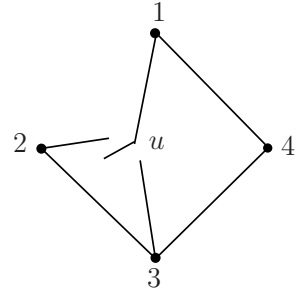}
\caption{Network Graph of MAS in Example \ref{e4.3.1}\label{fig4}}
\end{figure}

Consider a consensus problem of a multi-agent system with its network graph depicted in Fig. \ref{fig4}. There are 4 agents, $N=\{1,2,3,4\}$ with a common action set $S_i=\{1, 2\}$, $i=1,2,3,4$.  Assume all players can only  communicate with their neighbors. Additionally, there is a switch, denoted by $u$, which can link agent $1$ with $2$, or agent $1$ with $3$, or neither of them. The system objective function is
\begin{align}\label{4.3.1}
\phi(u, a)=2\dsum_{i\in N} {\bf 1}_{\{a_i=1\}}+\dsum_{(i,j)\in E(u)} \frac{{\bf 1}_{\{a_i=a_j\}}}{2} .
\end{align}
Define the state set $X=\{x_1,x_2,x_3\}$, where $x_1$ means the switch $u$ is open; $x_2$ means the switch $u$ is connected with node $3$; $x_3$ means the switch $u$ is connected with node $2$. Then there are $3$ states shown in Fig. \ref{fig5}.

\begin{figure}[!htbp]
\centering
\includegraphics[width=.65\linewidth]{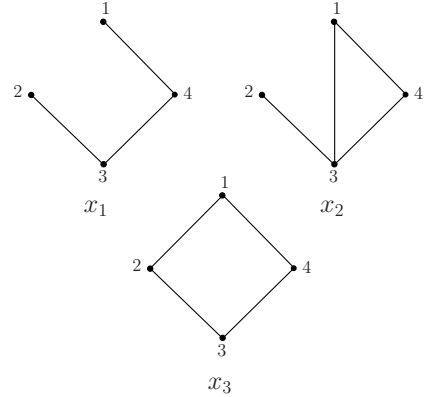}
\caption{Three states of the MAS\label{fig5}}
\end{figure}

Then we have
\begin{align}\label{4.3.2}
\phi(x_i, a)=V^{\phi(x_i, \cdot)} a,
\end{align}
where
\begin{align*}
  \begin{array}{ll}
  V^{\phi(x_1, \cdot)}&=[ 11  , 7  , 7  , 5  , 8  , 4  , 6  , 4  , 8  , 6  , 4  , 4  , 5  , 3  , 3  , 3],\\
   V^{\phi(x_2, \cdot)}&=[12  , 8  , 7  , 5  , 9  , 5  , 6  , 4  , 8  , 6  , 5  , 5  , 5  , 3  , 4  , 4],\\
  V^{\phi(x_3, \cdot)}&=[12  , 8  , 8  , 6  , 8  , 4  , 6  , 4  , 8  , 6  , 4  , 4  , 6  , 4  , 4  , 4].\\
  \end{array}
\end{align*}
Next, we design the SEP. Assume we use the SEP-2, then
according to (\ref{4.2.301}), we can design the state transition function as
\begin{align*}
V^{P(x_1,\cdot)}&=
\left[\begin{array}{cccccccccccccccc}
 \frac{1}{3}  & \frac{1}{3}  & \frac{1}{3}  & \frac{1}{3}  & \frac{1}{3}  & \frac{1}{3}  & \frac{1}{3}  & \frac{1}{3}  & \frac{1}{3}  & \frac{1}{3}  & \frac{1}{3}  & \frac{1}{3}  & \frac{1}{3}  & \frac{1}{3}  & \frac{1}{3}  & \frac{1}{3}\\
 \frac{1}{3}  & \frac{1}{3}  & \frac{1}{3}  & \frac{1}{3}  & \frac{1}{3}  & \frac{1}{3}  & \frac{1}{3}  & \frac{1}{3}  & \frac{1}{3}  & \frac{1}{3}  & \frac{1}{3}  & \frac{1}{3}  & \frac{1}{3}  & \frac{1}{3}  & \frac{1}{3}  & \frac{1}{3}\\
 \frac{1}{3}  & \frac{1}{3}  & \frac{1}{3}  & \frac{1}{3}  & \frac{1}{3}  & \frac{1}{3}  & \frac{1}{3}  & \frac{1}{3}  & \frac{1}{3}  & \frac{1}{3}  & \frac{1}{3}  & \frac{1}{3}  & \frac{1}{3}  & \frac{1}{3}  & \frac{1}{3}  & \frac{1}{3}
\end{array}\right],\\
V^{P(x_2,\cdot)}&=
\left[\begin{array}{cccccccccccccccc}
 0  & 0  & \frac{1}{3}  & \frac{1}{3}  & 0  & 0  & \frac{1}{3}  & \frac{1}{3}  & \frac{1}{3}  & \frac{1}{3}  & 0  & 0  & \frac{1}{3}  & \frac{1}{3}  & 0  & 0\\
 \frac{1}{2}  & \frac{1}{2}  & \frac{1}{3}  & \frac{1}{3}  & 1  & 1  & \frac{1}{3}  & \frac{1}{3}  & \frac{1}{3}  & \frac{1}{3}  & 1  & 1  & \frac{1}{3}  & \frac{1}{3}  & \frac{1}{2}  & \frac{1}{2}\\
 \frac{1}{2}  & \frac{1}{2}  & \frac{1}{3}  & \frac{1}{3}  & 0  & 0  & \frac{1}{3}  & \frac{1}{3}  & \frac{1}{3}  & \frac{1}{3}  & 0  & 0  & \frac{1}{3}  & \frac{1}{3}  & \frac{1}{2}  & \frac{1}{2}
\end{array}\right],
\end{align*}
\begin{align*}
V^{P(x_3,\cdot)}&=
\left[\begin{array}{cccccccccccccccc}
  0  & 0  & 0  & 0  & \frac{1}{3}  & \frac{1}{3}  & \frac{1}{3}  & \frac{1}{3}  & \frac{1}{3}  & \frac{1}{3}  & \frac{1}{3}  & \frac{1}{3}  & 0  & 0  & 0  & 0\\
 \frac{1}{2}  & \frac{1}{2}  & 0  & 0  & \frac{1}{3}  & \frac{1}{3}  & \frac{1}{3}  & \frac{1}{3}  & \frac{1}{3}  & \frac{1}{3}  & \frac{1}{3}  & \frac{1}{3}  & 0  & 0  & \frac{1}{2}  & \frac{1}{2}\\
 \frac{1}{2}  & \frac{1}{2}  & 1  & 1  & \frac{1}{3}  & \frac{1}{3}  & \frac{1}{3}  & \frac{1}{3}  & \frac{1}{3}  & \frac{1}{3}  & \frac{1}{3}  & \frac{1}{3}  & 1  & 1  & \frac{1}{2}  & \frac{1}{2}
 \end{array}\right].
 \end{align*}
Then we construct
$
\begin{bmatrix}\G_{x_i}\\ \B_{x_i}\end{bmatrix}, \; i=1,2, 3.
$

Note that $\B_{x_i}=\B$, which is the same as the one in Example \ref{e3.3.1}.
In addition, it is easy to calculate that
\begin{align*}
&\G_{x_1}=\begin{bmatrix}
I_2\ot {\bf 1}_4^T\ot I_2&~&~&\\
&{\bf 1}_2^T\ot I_4\ot{\bf 1}_2^T&~&~\\
 &~&{\bf 1}_2^T\ot I_8&\\
  &~&~&I_2\ot{\bf 1}_2^T\ot I_4
\end{bmatrix},\\
&\G_{x_2}=\begin{bmatrix}
I_2\ot {\bf 1}_2^T\ot I_4&~&~&\\
&{\bf 1}_2^T\ot I_4\ot{\bf 1}_2^T&~&~\\
 &~&I_{16}&\\
  &~&~&I_2\ot {\bf 1}_2^T\ot I_4
\end{bmatrix},\\
&\G_{x_3}=\begin{bmatrix}
I_4\ot {\bf 1}_2^T\ot I_2&~&~&\\
&I_8\ot{\bf 1}_2^T&~&~\\
 &~&{\bf 1}_2^T\ot I_8&\\
  &~&~&I_2\ot {\bf 1}_2^T\ot I_4
\end{bmatrix}.
\end{align*}
Using them, it is ready to verify that
\begin{align}\label{4.3.3}
{\bf 1}^T_4\ot V^{\phi(x_i, \cdot)}\in \SpanRow \begin{bmatrix}\G_{x_i}\\ \B_{x_i}\end{bmatrix}, \; i=1,2, 3.
\end{align}
It is easy to prove that (\ref{4.3.3}) is equivalent to (\ref{4.2.4}). According to Theorem \ref{t4.2.2}, there exists a state based potential game $G$ with a set of neighborhood-determinant utilities and the $\phi(x,a)$ in (\ref{4.3.1}) as its potential function.

A group of utilities functions are given bellow.
\begin{align}\label{4.3.5}
\begin{array}{l}
c_i(x,a)=2*{\bf 1}_{\{a_i=1\}}+\dsum_{j\in U^{x}(i)} {\bf 1}_{\{a_j=a_i\}},\\
 x\in X, a\in S, i=1,\cdots,4.
 \end{array}
\end{align}
Using the state transition function $P(x,a)$, we have the state dynamic equation
\begin{align}
  \label{4.3.5}
  x(t+1)=M_Px(t)a(t),
\end{align}
where
\[
M_P=[V^{P(x_1,\cdot)}, V^{P(x_2,\cdot)},V^{P(x_3,\cdot)}].
\]
Using better reply with inertia dynamics and let $\epsilon=0.1$, then we have the strategy dynamic equation
\begin{align}
  \label{4.3.5}
  a(t+1)=M_Fx(t+1)a(t),
\end{align}
where
\[
M_F=\begin{bmatrix}
1  & 0.9  & 0.9  & 0.81  & \dots  & 0  & 0\\
 0  & 0.1  & 0  & 0.09  & \dots  & 0  & 0\\
 0  & 0  & 0.1  & 0.09  & \dots  & 0  & 0\\
 0  & 0  & 0  & 0.01  & \dots  & 0  & 0\\
\vdots&&&&&&\\
 0  & 0  & 0  & 0  & \dots  & 0.09  & 0\\
 0  & 0  & 0  & 0  & \dots  & 0  & 0\\
 0  & 0  & 0  & 0  & \dots  & 0.01  & 0\\
 0  & 0  & 0  & 0  & \dots  & 0  & 1
\end{bmatrix}\in \M_{16\times 48}.
\]

Finally,  we can prove that $[a^*, x^*]$ is the only recurrent state equilibrium, where $a^*=(1,1,1,1)$ and $x^*\in \{x_2,x_3\}$. Then Theorem \ref{t4.1.4} assures that when the better reply with inertia dynamics is used, the action state pair converges almost surely to this  $[a^*, x^*]$. Note that  $a^*$ is unique and $x^*\in \{x_2,x_3\}$ is an invariant set. Fig. \ref{fig6} presents several simulations of the better reply inertia on the consensus problem of multi agent system. We can see that the action of all the four agents will reach the consensus $a^*=(1,1,1,1)$, which maximizes the objective function (\ref{4.3.1}) of the MAS.

 \begin{figure}[!htbp]
\centering
\begin{subfigure}[b]{.475\textwidth}
\includegraphics[width=\linewidth]{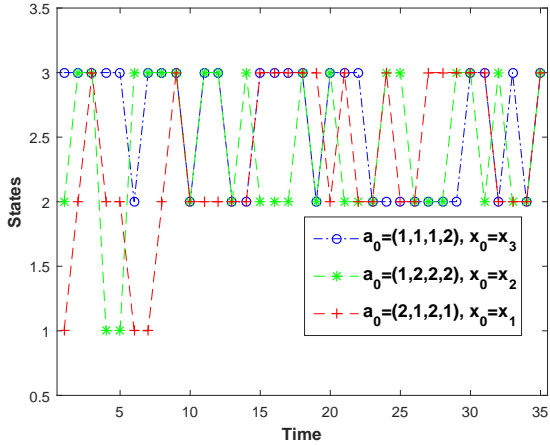}
\caption{State Dynamic\label{fig6.a}}
\end{subfigure}
\begin{subfigure}[b]{.475\textwidth}
\includegraphics[width=\linewidth]{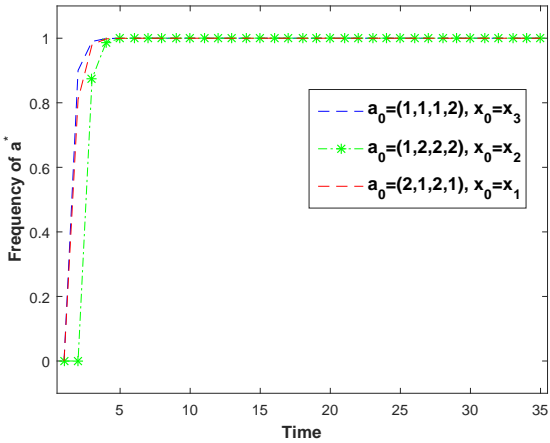}
\caption{Strategy Dynamic\label{fig6.a}}
\end{subfigure}
\caption{The dynamics of states and strategies in 35 times iterations initialled with $[(1,1,1,2),x_3]$, $[(1,2,2,2),x_2]$ and $[(2,1,2,1),x_1]$. Note that the better reply with inertia dynamic guarantees convergence to $a^*=(1,1,1,1)$.\label{fig6}}
\end{figure}
\end{exa}
\section{Conclusion}

This paper considers the problem of game theoretic control of multi-agent systems. Assume there is a system level objective function, say consensus, then we may design a local information based utility functions such that the multi-agent system becomes a potential game with the system level objective function as the potential function. Then individual agents can use their local information to reach an equilibrium, which maximizes the objective function. Two cases have been investigated. (i) fixed topology, and (ii) time-varying topology. Necessary and sufficient conditions have been obtained for both cases. Some examples are presented to illustrate the theoretical results.

There are several problems remaining for further study. For example,
\begin{enumerate}
\item The state transition process used in this paper for time-varying topology is designable. A challenge problem for further investigation is to relax this.

\item When the necessary and sufficient condition fails, can we design a near-potential game \cite{can13} to reach the same goal?
\end{enumerate}

\end{document}